\documentclass[11pt,a4paper]{amsart}
\usepackage[utf8]{inputenc}
\usepackage[T1]{fontenc}
\usepackage[english]{babel}
\usepackage{lmodern,microtype}
\usepackage{amsmath,amssymb,amsthm,mathtools}
\usepackage[margin=30mm]{geometry}
\usepackage{tikz}
\usetikzlibrary{patterns.meta}
\usepackage{placeins,needspace,booktabs,array}
\usepackage{verbatim}

\usepackage{color, xcolor}
\definecolor{darkgreen}{rgb}{0,0.5,0}
\usepackage[colorlinks, citecolor=darkgreen, backref]{hyperref}

\numberwithin{equation}{section}

\theoremstyle{plain}
\newtheorem{thm}{Theorem}[section]
\newtheorem{lm}[thm]{Lemma}
\newtheorem{prop}[thm]{Proposition}
\newtheorem{cor}[thm]{Corollary}
\theoremstyle{definition}

\newtheorem{ex}[thm]{Example}
\theoremstyle{remark}
\newtheorem{re}[thm]{Remark}
\newtheorem*{remark}{Remark}

\newcommand{\Z}{\mathbb{Z}}
\newcommand{\R}{\mathbb{R}}
\newcommand{\PP}{\mathbb{P}}
\newcommand{\EE}{\mathbb{E}}

\newcommand{\cP}{\mathcal{P}}
\newcommand{\bh}{\mathbf{h}}
\newcommand{\e}{\mathbf{e}}
\newcommand{\z}{\zeta(2)}
\newcommand{\Hh}{S_{\bh}}
\newcommand{\floor}[1]{\left\lfloor #1 \right\rfloor}

\date{November 26, 2026}

\title[Asymptotics for the rectangle and box partition functions]{Asymptotics for the rectangle and box partition functions}

\author{Krystian Gajdzica}
\address{Theoretical Computer Science Department \\ Faculty of Mathematics and Computer Science \\ Jagiellonian University \\ \L ojasiewicza 6 \\ 30-348 Krak\'ow \\ Poland}
\email{krystian.gajdzica@uj.edu.pl}

\author{Robin Visser}
\address{Charles University \\ Faculty of Mathematics and Physics \\ Department of Algebra \\ Sokolovsk\'a 83 \\ 186 75 Praha 8 \\ Czech Republic}
\email{robin.visser@matfyz.cuni.cz}

\author{Maciej Zakarczemny}
\address{Department of Applied Mathematics \\ Faculty of Computer Science and Mathematics \\ Cracow University of Technology \\ Warszawska 24 \\ 31-155 Krak\'ow \\ Poland}
\email{mzakarczemny@pk.edu.pl}

\author{B{\l}a{\.z}ej {\.Z}mija}
\address{Department of Mathematical Analysis and Applications \\ Faculty of Applied Mathematics \\ AGH University of Krakow \\ al. Mickiewicza 30 \\ 30-059 Krak\'ow \\ Poland}
\email{bzmija@agh.edu.pl}

\thanks{K. G. was supported by the National Science Center grant no. 2024/53/N/ST1/01538. R. V. was supported by {Charles University} programme PRIMUS/24/SCI/010 and {Charles University} programme UNCE/24/SCI/022.}



\makeatletter
\newcommand{\mathproofend}{\begingroup\def\math@qedhere{\quad\hbox{\qedsymbol}}\qedhere\endgroup}
\makeatother

\begin{document}
\begin{abstract}


In this paper, we study the box partition function $p(m_1,\ldots,m_d)$ that counts the multisets of boxes with positive integer
side lengths that admit an axis-parallel partition of the integer-sided box $[0,m_1]\times\cdots\times[0,m_d]$, where two partitions are indistinguishable if they consist of the same multiset of boxes up to coordinate permutations. For every $d\ge1$, we derive the asymptotics of $p(m_1,\ldots,m_{d-1},n)$ as $n$ approaches infinity. In particular, this resolves the recent conjecture by the first and the third authors concerning the asymptotic behaviour of rectangle partition function $p(m,n)$.\\
We also establish the logarithmic asymptotics of the box partition function when several side lengths grow as fixed positive powers of a common parameter.



\end{abstract}
\maketitle

\section{Introduction}

Let $p(n)$ be the number of \emph{partitions} of $n$, that is, the number of ways of writing $n$ as a sum of positive integers, without regard to ordering.  This function has played an important role in many areas of mathematics and other sciences. Throughout history, many prominent mathematicians like Andrews, Euler, Glaisher, Hardy and Ramanujan have contributed to its development.
In 1751, Euler obtained the following generating function for $p(n)$:
\begin{equation*}
    \sum_{n \ge 0} p(n) q^n = \prod_{n \ge 1} \tfrac{1}{1 - q^n} ,
\end{equation*}
and in 1918, Hardy and Ramanujan~\cite{HardyRamanujan} used the circle method to show the following explicit asymptotic formula for $p(n)$:
\begin{equation} \label{eq:HR}
    p(n) \sim \tfrac{1}{4 n \sqrt{3}} \exp  \bigg(  \pi \sqrt{\tfrac{2n}{3}} \bigg)
\end{equation}
as $n \to \infty$.  This was further improved by Rademacher \cite{Rademacher}, who proved an exact infinite series converging to $p(n)$.  We refer the reader to~\cite{Andrews, AndrewsEriksson, Apostol} and the references therein for further background on the partition function.

\subsection{Rectangle partitions}
Let $p(m,n)$ be the number of ways of partitioning a rectangle of size $m\times n$ into a multiset of rectangular blocks with positive integer sides. All sides are parallel to the coordinate axes. Blocks of sizes $a\times b$ and $b\times a$ are identified. Each multiset is counted once, regardless of how many arrangements realize it. In particular, $p(1,n)=p(n)$, where $p(n)$ is the usual partition function.

The values $p(m,n)$ appear in the OEIS as the sequence A360629, contributed by P. von Br\"{o}mssen in 2023~\cite{OEIS}. The values for $p(n,n)$, $p(2,n)$, and $p(3,n)$ are given in sequences A360630 to A360632~\cite{OEIS}.
This extension of the integer partition function was recently studied by the first three authors~\cite{GVZ}. They obtained an explicit asymptotic for $p(2, n)$ as $n \to \infty$, and conjectured~\cite[Conjecture~3.2]{GVZ} that
\begin{equation} \label{eq:conjp3n}
    p(3,n) \sim k_3 n^{\delta_3} \exp (\lambda_3 \sqrt{n} )
\end{equation}
for some real constants $k_3, \lambda_3 > 0$ and $\delta_3 \in \R$.
Their work also treats horizontally symmetric partitions and restrictions on the allowed block sizes. The first and third authors~\cite{GZsquare} subsequently determined the logarithmic asymptotics for $p(3,n)$ and $p(n,n)$, as $n \to \infty$.
More recently, they gave an explicit asymptotic for $\log p(m,n)$ as $n\to\infty$ \cite[Theorem~1.1]{GZfixed}, for any fixed $m\in\mathbb{Z}_+$.  In particular, they showed that
\begin{equation*}
    \log p(m, n) 
    = 
    \pi \sqrt{\tfrac{2m H_m n}{3}}
    + O(\log n),
\end{equation*}
as $n \to \infty$, where $H_m = \sum_{i=1}^m 1/i$ is the $m$-th harmonic number.  
In this paper, we strengthen this to an explicit asymptotic formula for $p(m,n)$ itself (Theorem~\ref{cor:rect}).  We obtain it as a special case of a more general result for higher-dimensional boxes.

\subsection{Higher-dimensional partitions}
In this paper, we shall focus mainly on extending the rectangular partition function and studying the more general $d$-dimensional box partition function $p(m_1, \dots, m_d)$ for some fixed dimension $d \ge 2$.  This counts the multisets of $d-$~dimensional boxes with positive integer side lengths that admit an axis-parallel partition of the integer-sided box $[0,m_1] \times \cdots \times [0,m_d]$, where two partitions are indistinguishable if they consist of the same multiset of boxes up to coordinate permutations.  Thus $d = 1$ gives the usual partition function $p(n)$ and $d = 2$ gives the rectangle partition function $p(m, n)$.

We first derive an asymptotic for boxes with a single-growing side. Let $d \ge 2$ and let $h_1 \le h_2 \le \dots \le h_{d-1}$ be fixed positive integers.  For brevity, we define $\mathbf{h} := (h_1, \dots, h_{d-1})$.  Our first theorem is an explicit asymptotic for $p(\mathbf{h}, n) = p(h_1, \dots, h_{d-1}, n)$ as $n$ goes to infinity.  We define the set
\begin{equation*}
    \mathcal{A}_\mathbf{h} := \left\{
        (a_1,\ldots,a_{d-1}) \in \mathbb{Z}_{+}^{d-1} :
         \begin{array}{l}
        a_1 \le \cdots \le a_{d-1},\\ 
        a_i \le h_i \text{ for all } i\in\{1,2,\ldots,d-1\}
        \end{array}
    \right\},
\end{equation*}
and define $A_\mathbf{h}, B_\mathbf{h}, C_\mathbf{h}$ as the real constants
\begin{align*}
    A_\mathbf{h} &:= \frac{\pi^2}{6}
    \sum_{
    \mathbf{a} 
    \in \mathcal{A}_\mathbf{h}}
    \frac{1}{a_1\cdots a_{d-1}}, \qquad 
    B_\mathbf{h} := \sum_{
    \mathbf{a} 
    \in \mathcal{A}_\mathbf{h}}
    \left(a_{d-1}-\frac{1}{2}\right), \quad \text{ and }\\[3pt]
    C_\mathbf{h} &:= \prod_{
    \mathbf{a} 
    \in \mathcal{A}_\mathbf{h}}
    \frac{(a_1\cdots a_{d-1})^{a_{d-1}-\tfrac12}(a_{d-1}-1)!}
    {\sqrt{2\pi}}.
\end{align*}
where the sums and the product range over all $\mathbf{a} = (a_1, \dots, a_{d-1})$ in $\mathcal{A}_{\mathbf{h}}$.  Then we can show the following theorem:

\begin{thm} \label{thm:fixed}
    For any fixed $d \ge 2$, and fixed positive integers $h_1 \le \dots \le h_{d-1}$, we have
    \begin{equation*}
        p(\mathbf{h}, n) = p(h_1, \dots, h_{d-1}, n) \sim \frac{C_\mathbf{h} A_\mathbf{h}^{B_{\mathbf{h}}/2+1/4}}{2\sqrt{\pi}}\,
    N_\mathbf{h}^{-B_{\mathbf{h}}/2-3/4}
    \exp\left(2\sqrt{A_\mathbf{h}N_\mathbf{h}}\right),
    \end{equation*}
    as $n \to \infty$, where $ N_\mathbf{h}=h_1\cdots h_{d-1} n$.
    
\end{thm}
\begin{remark}
When all fixed side lengths are equal, that is,
$h_1=\cdots=h_{d-1}=m$, the asymptotic formula in
Theorem~\ref{thm:fixed} takes the following compact form:
\begin{gather*}
p(m,\ldots,m,n)\\
\sim
\frac{(m!)^{(m-\frac12)\binom{m+d-2}{d-2}}
\pi^{(d-1)\binom{m+d-2}{d}}}
{2\sqrt{2^{\binom{m+d-2}{d-1}}m^{d-1}n}}
\cdot
\left(\frac{H_{m,d-1}}{6m^{d-1}n}\right)^{
\frac{2(d-1)\binom{m+d-2}{d}
+\binom{m+d-2}{d-1}+1}{4}}
\cdot
\exp\left(\pi\sqrt{\tfrac{2}{3}m^{d-1}H_{m,d-1}n}\right)
\end{gather*}
as $n\to\infty$, where
$
H_{m,d-1}:=
\sum_{j=1}^{m}(-1)^{j-1}\binom{m}{j}j^{1-d}.
$
In particular, $H_{m,1}=H_m$.
\end{remark}
For $d=2$ in Theorem~\ref{thm:fixed}, we have $\mathbf{h} = (m)$, $\mathcal{A}_{\mathbf{h}} = \{1, \dots, m \}$, $A_\mathbf{h} = \pi^2 H_m / 6$, $B_\mathbf{h} = m^2 / 2$, and $C_\mathbf{h} = (m!)^{m-\tfrac12} / (2 \pi)^{m/2}$, for some fixed $m \geq 1$.  This then gives an explicit asymptotic for the rectangular partition function $p(m,n)$ as $n \to \infty$:

\begin{thm}  \label{cor:rect}   
    For any fixed integer $m\ge1$, as $n\to\infty$, we have
    \begin{equation*}
    p(m,n)\sim
\frac{(m!)^{m-\frac12}\pi^{\binom{m}{2}}}
{2\sqrt{2^m mn}}
\left(\frac{H_m}{6mn}\right)^{\frac{m^2+1}{4}}
\exp\left(\pi\sqrt{\tfrac{2}{3}mH_mn}\right).
    \end{equation*}
\end{thm}

Theorem~\ref{cor:rect} recovers the classical Hardy--Ramanujan asymptotic \cite{HardyRamanujan} given in \eqref{eq:HR} for the usual partition function $p(n)$ when $m=1$ (with the convention $\binom{1}{2}=0$), and a theorem
of the first three authors~\cite[Theorem~1.5]{GVZ}
when $m=2$.  The case of $m = 3$ proves \eqref{eq:conjp3n} and gives the following new explicit asymptotic formula:
\begin{equation*}
    p(3, n) \sim \tfrac{121\pi^3}{7776}\sqrt{\tfrac{11}{3}}\;n^{-3}\exp\left(\pi\sqrt{\tfrac{11n}{3}}\right).
\end{equation*}
as $n \to \infty$, thereby strengthening the theorem of the first and third authors~\cite{GZsquare}, where only the logarithmic asymptotic for $p(3, n)$ was obtained.

We also mention that we've received recent communication from Tony Newton~\cite{Newton2026} who has independently obtained the asymptotic formula stated in Theorem~\ref{cor:rect}, by combining an application of Ingham's Tauberian theorem with a probabilistic packing argument.

\subsection{Several growing sides}
We can furthermore study the asymptotic behaviour for the $d$-dimensional partition function $p(m_1, \dots, m_d)$ where several sides grow simultaneously, as fixed powers of a common parameter $n$.
Let $r \ge 0$ and $g \ge 2$, let $h_1 \le h_2 \le \dots \le h_{r}$ be fixed positive integers, and let $0 < \alpha_1 \le \dots \le \alpha_g$ be fixed positive real numbers.
We define $S_{\mathbf{h}}$ as the real number
\begin{equation*}
    S_{\mathbf{h}} := \sum_{\substack{1\le a_1\le\cdots\le a_r\\
    a_i\le h_i\text{ for all }1\le i\le r}}
    \frac{1}{a_1\cdots a_r},
\end{equation*}
and let $D_{\mathbf{\alpha}}$ denote the Euclidean volume of the $(g-1)$-dimensional polytope:
\begin{equation*}
    \mathcal{D}_{\alpha} := 
        \left\{
        (y_1,\ldots,y_{g-1})\in\mathbb{R}^{g-1}:
        \begin{array}{l}
            0\le y_1\le\cdots\le y_g,\\
            y_i\le\alpha_i
            \text{ for all }i\in\{1,2,\ldots,g\}
        \end{array}
    \right\},
\end{equation*}
where
$
    y_g=\frac{1}{2}\sum_{i=1}^g\alpha_i
    -\sum_{i=1}^{g-1}y_i
$.  Then we have the following logarithmic asymptotic for $p(h_1,\ldots,h_r,\floor{n^{\alpha_1}}, \dots, \floor{n^{\alpha_g}})$:

\begin{thm} \label{thm:growing} 
    For any fixed $r \ge 0$, $g \ge 2$, 
    $h_1 \le \dots \le h_r$, and $0 < \alpha_1 \le \dots \le \alpha_g$ as above,
    we have
    \begin{equation*}
    \log p \big( h_1,\ldots,h_r,
    \floor{n^{\alpha_1}}, \dots, \floor{n^{\alpha_g}}
    \big)
    \sim
    \pi\sqrt{\frac{2 S_\mathbf{h} D_{\boldsymbol{\alpha}}N}{3}}\,
    (\log n)^{(g-1)/2}.
    \end{equation*}
    as $n \to \infty$, where $N = h_1 \cdots h_r \cdot \floor{n^{\alpha_1}} \cdots \floor{n^{\alpha_g}}$.
\end{thm}

The above theorem
can be specialised to many particular cases.  For example, for any $d \ge 2$, we can take $r = 0$ and $\alpha_1 = \dots = \alpha_g = 1$ to obtain the following logarithmic asymptotics for the growing hypercubes in dimensions two, three, and four:
\begin{align*}
 \log p(n,n)\sim\tfrac{\pi}{\sqrt3}n\sqrt{\log n},\
 \log p(n,n,n)\sim\tfrac{\pi}{2\sqrt3}n^{\frac32}\log n,\
 \log p(n,n,n,n)\sim\tfrac{\pi}{3\sqrt6}n^2(\log n)^{\frac32}.
\end{align*}
as $n \to \infty$.  Whilst the first of these was obtained in~\cite{GZsquare}, this gives new results in all higher dimensions.  In general, we obtain the following neat corollary for the $d$-dimensional hypercube:

\begin{cor} \label{cor:cubes}
    For each fixed $d \ge 2$, we have
    \begin{equation*}
        \log p (n, n, \dots, n) \sim \kappa_d \,n^{d/2}(\log n)^{(d-1)/2}
    \end{equation*}
    as $n \to \infty$ (with all $d$ arguments in the $d$-dimensional partition function equal to $n$), and where
    \begin{equation*}
        \kappa_d := \pi\sqrt{
        \tfrac{2}{3d!(d-1)!}
        \sum_{j=0}^{\lfloor d/2\rfloor}
        (-1)^j\binom{d}{j}
        \left(\tfrac{d}{2}-j\right)^{d-1}}.
    \end{equation*}
    
\end{cor}
For rectangles with two growing sides, the first and third authors conjectured~\cite[Remark~5.9]{GZsquare} that if $h\le w$ and $h\to\infty$, then
\begin{equation*}
    \Bigl(\tfrac{\pi}{\sqrt3}-o(1)\Bigr)\sqrt{hw\log h}\;\le\;\log p(h,w)\;\le\;\Bigl(\tfrac{\pi}{\sqrt6}+o(1)\Bigr)\sqrt{hw\log(hw)},
\end{equation*}
and that $\log p(h,w)\sim\frac{\pi}{\sqrt3}\sqrt{hw\log h}$ if, moreover, $\log w\sim\log h$. The following corollary confirms these conjectures when $h$ and $w$ grow polynomially in a common parameter, and identifies the correct constant.


\begin{cor} \label{cor:special}
    Fix some $0 < s \le u$, and let $h = \floor{n^s}$ and $w = \floor{n^u}$. Then
    \begin{equation}
        \log p (h, w) 
        \sim 
        \begin{cases}
            \tfrac{\pi}{\sqrt6}\sqrt{hw\log(hw)}, & \text{if } u \le 3s, \\[3mm]
            \pi\sqrt{\tfrac{2}{3}}\sqrt{hw\log h}, & \text{if } u \ge 3s ,
        \end{cases}
    \end{equation}
    as $n \to \infty$.
\end{cor}

In particular, the conjectured upper bound is attained precisely when $u\le3s$, that is, when $w\le h^{3+o(1)}$, and the case $s=u$ gives the conjectured asymptotic when $\log w\sim\log h$. The second case is also consistent with Theorem~\ref{cor:rect}, which gives $\log p(h,w)\sim\pi\sqrt{2hH_hw/3}$ for fixed $h$, since $H_h\sim\log h$ as $h\to\infty$.

\subsection{Methods} \label{sec:methods}
We provide here a brief overview of the methods used in our paper.  
For Theorem~\ref{thm:fixed}, put $r=d-1$ and let $U_\bh(N)$ be the number of multisets of total volume $N$ of boxes with sorted side lengths $(a_1,\ldots,a_r,c)$, where $a\in\mathcal{A}_\bh$ and $c\ge a_r$. Neither the condition $c\le n$ nor any tiling condition is imposed. We note that its generating function is
\begin{equation} \label{eq:Uh-gf}
    \sum_{N \ge 0} U_\bh(N) \, q^N = \prod_{a\in\mathcal{A}_\bh}\ \prod_{c\ge a_r}\frac{1}{1-q^{a_1\cdots a_r c}},
\end{equation}
and every partition of the $h_1\times\cdots\times h_r\times n$ box is such a multiset, so the upper bound $p(\bh,n)\le U_\bh(N)$ with $N= h_1 \cdots h_{d-1} n$ as in Theorem~\ref{thm:fixed}, is immediate. The asymptotic for $U_\bh(N)$, which is the right-hand side of Theorem~\ref{thm:fixed}, follows from Ingham's Tauberian theorem~\cite{Ingham}, in the form given by Bringmann, Jennings-Shaffer and Mahlburg~\cite{BJM}, applied to the product~\eqref{eq:Uh-gf} (Proposition~\ref{prop:coef}). The substance of Theorem~\ref{thm:fixed} is thus the lower bound $p(\bh,n)\ge(1-o(1))\,U_\bh(N)$, which we prove in Theorem~\ref{thm:ratio}. 


To obtain this lower bound, we use the Boltzmann probabilistic model~\cite{DFLS04,Vershik96}, in which the multiplicity of each type of box is an independent geometric random variable. Conditioned on the total volume being $N$, this is the uniform distribution on the multisets counted by $U_\bh(N)$, as in Fristedt's work on ordinary partitions~\cite{Fristedt}. We show that, with probability tending to one, such a multiset has the following properties: it contains no very long pieces; for each cross-section, the total length of the pieces with that cross-section is close to its expected value; and it contains many short rods and unit cubes (Lemma~\ref{lm:typical}). We then give an explicit construction showing that every multiset with these properties tiles the $h_1\times\cdots\times h_r\times n$ box (Section~\ref{sec:fixed-proof}). Its main ingredient is a random subdivision of the cross-section in which each piece cross-section appears, on average, with the required frequency (Lemma~\ref{lm:patterns}).

For boxes with several growing sides, the upper bound in Theorem~\ref{thm:growing} again follows easily, from a saddle-point bound on the number of multisets of individually fitting pieces. For the lower bound, we restrict ourselves to pieces whose growing side lengths are small relative to the corresponding side lengths of the box. We then round their sides to a geometric grid and pack them in layers.

\medskip
This paper is organised as follows. Section~\ref{sec:prelim} contains preliminaries, including the asymptotic for $U_\bh(N)$ and the Boltzmann model. Section~\ref{sec:packing} collects packing lemmas and constructs the cross-sectional subdivisions. Section~\ref{sec:fixed-proof} proves Theorem~\ref{thm:fixed} and Theorem~\ref{cor:rect}, and Section~\ref{sec:growing} proves Theorem~\ref{thm:growing} and Corollaries~\ref{cor:cubes} and~\ref{cor:special}. 
Finally, Section~\ref{sec:problems} lists some open problems.


\section{Preliminaries}\label{sec:prelim}

In this section we fix the notation for Sections~\ref{sec:prelim}--\ref{sec:fixed-proof}, determine the asymptotic behaviour of $U_\bh(N)$ (Proposition~\ref{prop:coef}), and introduce the Boltzmann probabilistic model. The main result of the section is Lemma~\ref{lm:typical}, which lists four properties that a uniformly random multiset counted by $U_\bh(N)$ has with probability tending to one. These are the properties used in Section~\ref{sec:fixed-proof} to show that such a multiset tiles the box.

\subsection{Pieces and types}
Throughout Sections~\ref{sec:prelim}--\ref{sec:fixed-proof} we fix $d\ge2$ and $\bh=(h_1,\ldots,h_{d-1})$ as in Theorem~\ref{thm:fixed}. To lighten notation we write $r=d-1$, so that $\bh=(h_1,\ldots,h_r)$, and we abbreviate $A=A_\bh$, $B=B_\bh$ and $C=C_\bh$. For $n\ge1$ we call
\[
[0,h_1]\times\cdots\times[0,h_r]\times[0,n]
\]
the \emph{host box}. This is the box to be tiled: $p(\bh,n)$ is the number of multisets of boxes that tile the host box. To avoid confusion with the host box, we call the boxes in such a multiset \emph{pieces}. We call $[0,h_1]\times\cdots\times[0,h_r]$ the \emph{cross-section} of the host box, and its last coordinate the \emph{axis}. The cross-section has volume $V=h_1\cdots h_r$, so the host box has volume $N=Vn$.

Since pieces are identified up to permutations of their side lengths, only their sorted side lengths matter. 
We start by noting the following easy observation concerning when one box fits inside another.

\begin{lm}\label{lm:fit}
Let $s_1\le\cdots\le s_d$ and $t_1\le\cdots\le t_d$ be positive integers. An $s_1\times\cdots\times s_d$ box admits an axis-parallel placement inside a $t_1\times\cdots\times t_d$ box if and only if $s_i\le t_i$ for every $i$.
\end{lm}

\begin{proof}
It's clear that $s_i\le t_i$ for every $i$ is a sufficient condition. To show it's a necessary condition, we note that an axis-parallel placement matches each $s_i$ with a distinct $t_{\pi(i)}\ge s_i$. For each $i$, the values $s_i,\ldots,s_d$ are matched with $d-i+1$ distinct sides of length at least $s_i$, so at least $d-i+1$ of the $t$'s are at least $s_i$, which gives $t_i\ge s_i$.
\end{proof}

Let $n\ge h_r$, so that $h_1\le\cdots\le h_r\le n$ are the sorted side lengths of the host box. By Lemma~\ref{lm:fit}, a piece with sorted side lengths $s_1\le\cdots\le s_d$ fits inside the host box if and only if $(s_1,\ldots,s_r)\in\mathcal{A}_{\bh}$ and $s_d\le n$. We therefore describe a piece by its \emph{type} $(a,c)$: its \emph{footprint} $a=(s_1,\ldots,s_r)$ lists its $r$ shortest sides, and its \emph{length} $c=s_d$ is its longest side. For $a\in\mathcal{A}_{\bh}$ we write
\[
v_a=a_1\cdots a_r\qquad\text{and}\qquad b_a=a_r,
\]
so that $\Hh=\sum_{a\in\mathcal{A}_{\bh}}v_a^{-1}$ and $A=\z\Hh$. A piece of type $(a,c)$ has volume $v_ac$, and the pieces that fit inside the host box are those of type $(a,c)$ with $a\in\mathcal{A}_{\bh}$ and $b_a\le c\le n$.

The words \emph{footprint} and \emph{length} refer to the sorted side lengths of a piece, not to its position: in a partition of the host box, a piece may be placed in any orientation. In the partitions that we construct in Section~\ref{sec:fixed-proof}, however, every piece of type $(a,c)$ lies with its footprint in the cross-section and its length along the axis. The \emph{unit footprint} is $\e=(1,\ldots,1)\in\mathcal{A}_{\bh}$. The pieces of type $(\e,c)$ are the \emph{rods} $1\times\cdots\times1\times c$, and the rod of type $(\e,1)$ is the \emph{unit cube}.

\medskip
Recall from Section~\ref{sec:methods} that $U_\bh(N)$ is the number of multisets of total volume $N$ of pieces of types $(a,c)$ with $a\in\mathcal{A}_{\bh}$ and $c\ge b_a$. Compared with $p(\bh,n)$, we drop both the upper bound $c\le n$ on the lengths and the requirement that the multiset tiles the host box. The generating function~\eqref{eq:Uh-gf} of $U_\bh(N)$ is
\begin{equation}\label{eq:Fh}
\Phi_{\bh}(q)=\prod_{a\in\mathcal{A}_{\bh}}\prod_{c=b_a}^{\infty}\frac1{1-q^{v_ac}}=\sum_{N\ge0}U_{\bh}(N)\,q^N.
\end{equation}
Every partition of the host box is one of the multisets counted by $U_\bh(Vn)$, so
\begin{equation}\label{eq:trivial-upper}
p(\bh,n)\le U_\bh(Vn)\qquad(n\ge h_r).
\end{equation}

\subsection{The generating product}
To apply a Tauberian theorem to $\Phi_\bh$, we need the behaviour of $\Phi_\bh(q)$ as $q\to1$. We write $q=e^{-z}$ and let $z\to0$ in a sector of the right half-plane. Let $P(q)=\prod_{j\ge1}(1-q^j)^{-1}$ be the generating function of $p(n)$, and for $\Delta>0$ let $\Sigma_\Delta=\{z\in\mathbb C:\ \operatorname{Re}z>0,\ |\operatorname{Im}z|\le\Delta\operatorname{Re}z\}$. Powers of $z$ are taken with the principal branch.
We prove an auxiliary result concerning the asymptotics of considered generating functions that will be used throughout this section.

\begin{lm}\label{lm:product}
Fix $\Delta>0$. Uniformly for $z\in \Sigma_\Delta$ with $z\to0$:
\begin{enumerate}
\item[(a)] $P(e^{-z})=\sqrt{z/(2\pi)}\,e^{\z/z}(1+O(|z|))$;
\item[(b)] for fixed real $v>0$ and integer $b\ge1$,
\[
\prod_{c\ge b}\frac1{1-e^{-vcz}}=\frac{v^{b-\tfrac12}(b-1)!}{\sqrt{2\pi}}\,z^{b-\tfrac12}\,e^{\z/(vz)}\bigl(1+O(|z|)\bigr);
\]
\item[(c)] $\Phi_{\bh}(e^{-z})=C\,z^{B}\,e^{A/z}(1+O(|z|))$.
\end{enumerate}
\end{lm}
Here powers of $z$ use the principal branch.
\begin{proof}
The modular transformation of the Dedekind eta function (see, e.g.,~\cite[Chapter~3]{ApostolModular}) gives
\[
P(e^{-z})=\sqrt{\frac z{2\pi}}\exp\left(\frac{\pi^2}{6z}-\frac z{24}\right)P\bigl(e^{-4\pi^2/z}\bigr)\qquad(\operatorname{Re}z>0).
\]
For $z\in \Sigma_\Delta$ we have $\operatorname{Re}(1/z)=\operatorname{Re}z/|z|^2\ge(1+\Delta^2)^{-\tfrac12}|z|^{-1}$, so $P(e^{-4\pi^2/z})=1+O(e^{-c/|z|})$ for some $c=c(\Delta)>0$. This proves (a). For (b), write the product as $P(e^{-vz})\prod_{c<b}(1-e^{-vcz})$ and use $1-e^{-w}=w(1+O(|w|))$ for the finitely many factors. Since $\prod_{c<b}vcz=v^{b-1}(b-1)!\,z^{b-1}$, (b) follows from (a) applied at $vz\in \Sigma_\Delta$. Multiplying (b) over $a\in\mathcal{A}_{\bh}$, with $v=v_a$ and $b=b_a$, gives (c).
\end{proof}

\subsection{Coefficient asymptotics}
The first application of Lemma \ref{lm:product} exposes in the proof of the asymptotic behaviour of the function $U_{\bh}(N)$. However, the main tool to deal with that issue is Ingham's Tauberian theorem due to Bringmann et al. \cite{BJM}.
\begin{prop}\label{prop:coef}
As $N\to\infty$ through all integers,
\[
U_{\bh}(N)\sim\frac{C}{2\sqrt\pi}\,A^{\frac B2+\frac14}\,N^{-\frac B2-\frac34}\,e^{2\sqrt{AN}}.
\]
\end{prop}

\begin{proof}
We apply Ingham's Tauberian theorem~\cite{Ingham} in the form of~\cite[Theorem~1.1]{BJM}. Its hypotheses are as follows. The coefficients $U_{\bh}(N)$ are nonnegative and weakly increasing, because $\Phi_{\bh}$ contains the factor $(1-q)^{-1}$ of the unit cube. Next, $\Phi_{\bh}(e^{-t})\sim Ct^Be^{A/t}$ as $t\to0^+$. Finally, $|\Phi_{\bh}(e^{-z})|\ll|z|^Be^{A/|z|}$ in every sector $\Sigma_\Delta$. The last two follow from Lemma~\ref{lm:product}(c), since $|e^{A/z}|=e^{A\operatorname{Re}(1/z)}\le e^{A/|z|}$. The theorem gives the stated asymptotic.
\end{proof}

\begin{re}\label{re:meinardus}
Proposition~\ref{prop:coef} can also be derived from Meinardus' theorem~\cite{Meinardus} (see~\cite[Theorem~6.2]{Andrews} or~\cite[Theorem~1]{GSE}). The associated Dirichlet series is

\[
\sum_{a\in\mathcal{A}_{\bh}}\sum_{c\ge b_a}(v_ac)^{-s}=\zeta(s)\sum_{a\in\mathcal{A}_{\bh}}v_a^{-s}-\sum_{a\in\mathcal{A}_{\bh}}\sum_{c<b_a}(v_ac)^{-s}.
\]
It has a simple pole at $s=1$ with residue $\Hh$, takes the value $-B$ at $s=0$, and has derivative $\log C$ there. A self-contained probabilistic proof is also possible: it uses a local limit theorem for the total volume in the model of Section~\ref{sec:boltzmann}, as in~\cite{Fristedt,GSE}.
\end{re}

\begin{re}\label{re:elementary-upper}
For $\bh=(m)$, the elementary bound $U_{\bh}(N)\le e^{Nt}\Phi_{\bh}(e^{-t})$, together with $\log\prod_{b\ge1}(1-e^{-abt})^{-1}\le\pi^2/(6at)$ and the choice $t=\pi\sqrt{H_m/(6mn)}$, gives
\[
p(m,n)\le\exp\left(\pi\sqrt{\frac{2mH_mn}{3}}\right)\qquad(n\ge m),
\]
which is the upper bound of~\cite{GZfixed}.
\end{re}

\subsection{The Boltzmann model}\label{sec:boltzmann}
In this subsection, we introduce the probabilistic model that plays the crucial role in the proofs of the main results.
To prove that almost all multisets counted by $U_\bh(N)$ tile the host box, we need to understand what a uniformly random such multiset looks like. This is easier in a model in which the numbers of pieces of different types are independent. We use the Boltzmann model, as Fristedt~\cite{Fristedt} did for ordinary partitions (see also~\cite{DFLS04,Vershik96}).

\medskip

Fix $\bh$ and $t>0$. 
Let $X_{a,c}$, for $a\in\mathcal{A}_{\bh}$ and $c\ge b_a$, be independent random variables with
\[
\PP_t(X_{a,c}=k)=\bigl(1-e^{-tv_ac}\bigr)e^{-tv_ack}\qquad(k=0,1,2,\ldots).
\]
We read $X_{a,c}$ as the multiplicity of the type $(a,c)$. Put
\[
T=\sum_{a,c}v_ac\,X_{a,c},\qquad W_a=\sum_{c\ge b_a}c\,X_{a,c},
\]
so that $T$ is the total volume and $W_a$ the total length of the pieces with footprint $a$. Since $\EE_tT<\infty$, almost surely only finitely many $X_{a,c}$ are nonzero. A multiset $\mu$ of volume $N$ has probability $e^{-tN}/\Phi_{\bh}(e^{-t})$, so
\begin{equation}\label{eq:conditioning}
\PP_t(T=N)=\frac{U_{\bh}(N)e^{-tN}}{\Phi_{\bh}(e^{-t})},
\end{equation}
and conditionally on $T=N$ the multiset is uniform among the $U_{\bh}(N)$ multisets of volume $N$. From now on we take
\[
t=t_N=\sqrt{A/N}.
\]
Then Proposition~\ref{prop:coef}, Lemma~\ref{lm:product}(c) and~\eqref{eq:conditioning} give
\[
\PP_t(T=N)\sim\frac{C}{2\sqrt\pi}A^{\frac B2+\frac14}N^{-\frac B2-\frac34}\cdot\frac{t^{-B}}{C}
=\frac{t^{3/2}}{2\sqrt{\pi A}}.
\]
Consequently, for all large $N$ and every event $E$,
\begin{equation}\label{eq:conditional-loss}
\PP_t(E\mid T=N)\le\frac{\PP_t(E)}{\PP_t(T=N)}\le3\sqrt{\pi A}\;t^{-3/2}\,\PP_t(E).
\end{equation}

\begin{lm}\label{lm:typical}
Fix $\bh$, an integer $R\ge1$ and $\varepsilon>0$. Let $N\to\infty$, $t=t_N$, and put
\[
L=\left\lceil\frac{6\log(1/t)}{t}\right\rceil,\qquad w=\lfloor t^{-\tfrac13}\rfloor,\qquad x_a=\frac{V}{\Hh v_a^2}\quad(a\in\mathcal{A}_{\bh}).
\]
Then each of the following events has $\PP_t(\,\cdot\mid T=N)\to1$:
\begin{enumerate}
\item[(i)] $X_{a,c}=0$ for every $a$ and every $c>L$;
\item[(ii)] $W_a\le(x_a+\varepsilon)N/V$ for every $a\in\mathcal{A}_{\bh}$;
\item[(iii)] $Y:=\sum_{c=2}^{w}c\,X_{\e,c}>RL$;
\item[(iv)] $X_{\e,1}\ge Rw$.
\end{enumerate}
\end{lm}

\begin{proof}
By~\eqref{eq:conditional-loss}, for (i)--(iii) it suffices to show that the unconditional probabilities of the complements are $o(t^{3/2})$.

(i) For small $t$, $1-e^{-t}\ge t/2$, so
\[
\sum_{a}\sum_{c>L}\PP_t(X_{a,c}\ge1)=\sum_a\sum_{c>L}e^{-tv_ac}\le|\mathcal{A}_{\bh}|\frac{e^{-tL}}{1-e^{-t}}\le2|\mathcal{A}_{\bh}|\,t^{5}.
\]

(ii) Fix $a$ and $0<s<v_a$, and write $v=v_a$, $b=b_a$. By independence,
\[
\EE_te^{stW_a}=\prod_{c\ge b}\frac{1-e^{-tvc}}{1-e^{-t(v-s)c}}.
\]
Lemma~\ref{lm:product}(b) at the real points $(v-s)t$ and $vt$ gives
\[
\log\EE_te^{stW_a}=\frac{\z}{t}\left(\frac1{v-s}-\frac1v\right)+O_{a,s}(1).
\]
Since $t^2=A/N$ and $x_aA/V=\z/v^2$, we have $st(x_a+\varepsilon)N/V=(s\z/v^2+s\varepsilon A/V)/t$. By Markov's inequality,
\[
\PP_t\bigl(W_a>(x_a+\varepsilon)N/V\bigr)\le\exp\left(\frac1t\left[\z\Bigl(\frac1{v-s}-\frac1v-\frac s{v^2}\Bigr)-\frac{s\varepsilon A}{V}\right]+O_{a,s}(1)\right).
\]
The term in round brackets equals $s^2/(v^2(v-s))$. Choosing $s$ small enough makes the square bracket at most $-s\varepsilon A/(2V)$, so the probability is $O(e^{-c/t})$ for some $c>0$.

(iii) We have $\EE_te^{-tcX_{\e,c}}=(1+e^{-tc})^{-1}$. Since $tw\le t^{2/3}\to0$, eventually $1+e^{-tc}\ge1+e^{-tw}>e^{\tfrac12}$ for all $c\le w$. Hence
\[
\PP_t(Y\le RL)\le e^{tRL}\,\EE_te^{-tY}\le e^{R+6R\log(1/t)}e^{-(w-1)/2},
\]
which is $o(t^{3/2})$ because $w\asymp t^{-\tfrac13}$.

(iv) Deleting $j$ unit cubes is a bijection from the multisets of volume $N$ containing at least $j$ unit cubes onto the multisets of volume $N-j$. Hence $\PP_t(X_{\e,1}\ge j\mid T=N)=U_{\bh}(N-j)/U_{\bh}(N)$. For $j=Rw=O(N^{1/6})$, Proposition~\ref{prop:coef} shows that this ratio is
\[
(1+o(1))\exp\bigl(2\sqrt A(\sqrt{N-j}-\sqrt N)\bigr)\longrightarrow1.\qedhere
\]
\end{proof}

\section{Packing lemmas and cross-sectional patterns}\label{sec:packing}

In this section, we develop packing arguments that are heavily used to deduce the exact asymptotics of considered partition statistics.
The following lemma is taken from bin packing algorithms developed by Johnson~\cite{Johnson74}:

\begin{lm}\label{lm:nextfit}
Let $\ell_1,\ldots,\ell_k\ge0$ and $L\ge1$ be integers, and let $s_1,\ldots,s_q$ be positive integers, each at most $L$, with $\sum_is_i\le\sum_j\ell_j-kL$. Then the $s_i$ can be distributed among $k$ bins of capacities $\ell_1,\ldots,\ell_k$, so that the items in each bin can be laid end to end from the start of the bin.
\end{lm}
\begin{proof}
Place the items in order, moving to the next bin whenever the next item does not fit. Every bin that is left behind has unused capacity less than $L$. If some item remained after the last bin, the placed items would have total more than $\sum_j\ell_j-kL$, a contradiction.
\end{proof}

\begin{lm}\label{lm:exact}
Let $L,w,k\in\Z^+$ and let $\ell_1,\ldots,\ell_k\ge0$ be integers. Suppose a multiset of positive integers, each at most $L$, has sum $\ell_1+\cdots+\ell_k$, that its elements in $\{2,\ldots,w\}$ have sum greater than $kL$, and that at least $kw$ of its elements equal $1$. Then it can be split into $k$ groups with sums $\ell_1,\ldots,\ell_k$.
\end{lm}
\begin{proof}
First place the elements larger than $w$ by next fit. If one were left over, every bin would have unused capacity less than $L$, so the total unused capacity would be less than $kL$. But that capacity equals the total of the unplaced elements, which exceeds $kL$ because of the elements in $\{2,\ldots,w\}$. Next place the elements in $\{2,\ldots,w\}$ by next fit into the remaining capacities. A failure would leave less than $kw$ free in total, but at least $kw$ ones remain. Finally, the ones fill the remaining capacity exactly.
\end{proof}

For the sake of convenience, let us introduce some auxiliary notation regarding tiling.

A \emph{pattern} for $\bh$ is a tiling of $[0,h_1]\times\cdots\times[0,h_r]$ by boxes $\prod_{i}[z_i,z_i+a_i]$ with integer $z_i$ and $a\in\mathcal{A}_{\bh}$; the box has side $a_i$ along axis $i$. We call these boxes the \emph{slots} of the pattern, and write $\nu_a(P)$ for the number of slots of footprint $a$ in a pattern $P$. A uniformly chosen unit cell of the cross-section lies in a slot of footprint $a$ with probability $v_a\nu_a(P)/V$.

\begin{lm}\label{lm:interval}
Let $1\le u\le M$ be integers and $q_1\ge\cdots\ge q_u\ge0$ with $\sum_iq_i=1$. There is a random subdivision of $[0,M]$ into integer intervals such that, for a uniformly random unit cell chosen independently, the interval containing it has length $i$ with probability $q_i$. If $u=1$ or $q_1>q_2$, the subdivision into unit intervals has positive probability.
\end{lm}
\begin{proof}
First let $q$ be uniform on $\{1,\ldots,j\}$ for some $j\le M$. If $M=j$, take a uniform permutation of $\{1,\ldots,M\}$ and lay its cycles out as consecutive intervals of lengths equal to the cycle lengths. The expected number of cycles of length $i$ is $\binom Mi(i-1)!(M-i)!/M!=1/i$ (see also~\cite{SheppLloyd}), so a uniform cell lies in an interval of length $i$ with probability $i\cdot\frac1i\cdot\frac1M=\frac1M$. If $M>j$, condition the permutation on having at least two cycles; the cell's interval length becomes uniform on $\{1,\ldots,M-1\}$, since the single $M$-cycle has probability $1/M$. Then subdivide each interval of length greater than $j$ recursively, using the construction for that length. By induction on $M$, a cell in such an interval sees each length in $\{1,\ldots,j\}$ with probability $1/j$, so the overall probability of each $i\le j$ is
\[
\frac1{M-1}+\frac{M-1-j}{M-1}\cdot\frac1j=\frac1j.
\]
For general $q$, use the construction for $\{1,\ldots,i\}$ with probability $i(q_i-q_{i+1})$, where $q_{u+1}=0$. These weights are nonnegative and sum to $\sum_iq_i=1$, and the probability of length $j$ telescopes to $\sum_{i\ge j}(q_i-q_{i+1})=q_j$. The case $i=1$ is the unit subdivision, with weight $q_1-q_2$.
\end{proof}

\begin{lm}\label{lm:patterns}
There is a random pattern $\mathsf P$ for $\bh$, taking finitely many values, such that
\[
\EE\,\nu_a(\mathsf P)=x_a=\frac{V}{\Hh v_a^2}\quad(a\in\mathcal{A}_{\bh})\qquad \text{and}\qquad\PP(\mathsf P\text{ consists of unit cells})>0.
\]
\end{lm}
\begin{proof}
Put $g_1(j)=\mathbf1_{j\le h_1}/j$ and, for $2\le i\le r$,
\[
g_i(j)=\frac{\mathbf1_{j\le h_i}}{j}\sum_{k=1}^{j}g_{i-1}(k).
\]
By induction, $g_i(j)$ is the sum of $1/(a_1\cdots a_i)$ over $1\le a_1\le\cdots\le a_i=j$ with $a_l\le h_l$; in particular $\sum_jg_r(j)=\Hh$. Each $g_i$ is weakly decreasing, because averages of initial segments of a weakly decreasing sequence decrease, and $g_i(1)>g_i(2)$.

Cut axis $r$ by Lemma~\ref{lm:interval} with $q_j=g_r(j)/\Hh$. Then cut each resulting slab of width $j$ independently along axis $r-1$, with $q_k=g_{r-1}(k)/(jg_r(j))$ for $k\le j$; these numbers sum to one, because $\sum_{k\le j}g_{r-1}(k)=jg_r(j)$ for $j\le h_r$. Continue in the same way down to axis $1$. Every cut uses a weakly decreasing law supported on at most the length of the axis, so Lemma~\ref{lm:interval} applies. The resulting slots have sides $a_1\le\cdots\le a_r$ along the axes $1,\ldots,r$, with $a_i\le h_i$.

The coordinates of a uniform unit cell are independent of each other and of the cuts. Hence the probability that the cell lies in a slot of footprint $a$ telescopes:
\[
\frac{g_r(a_r)}{\Hh}\prod_{i=1}^{r-1}\frac{g_i(a_i)}{a_{i+1}g_{i+1}(a_{i+1})}=\frac{g_1(a_1)}{\Hh\,a_2\cdots a_r}=\frac1{\Hh v_a}.
\]
Since this probability equals $v_a\EE\nu_a(\mathsf P)/V$, we obtain $\EE\nu_a(\mathsf P)=x_a$. With positive probability the first cut uses only unit intervals, and then all later cuts are forced to do the same.
\end{proof}

\begin{figure}[htbp]
\centering
\begin{tikzpicture}[x=1cm,y=1cm,
  zone/.style={draw,line width=.55pt},
  hatched/.style={pattern={Lines[distance=5pt,angle=45,line width=.3pt]},pattern color=black!40},
  lab/.style={font=\small,fill=white,fill opacity=.9,text opacity=1,inner sep=2pt}]
\def\u{0.62}
\def\xA{0}
\def\xF{13}
\pgfmathsetmacro{\xB}{13*5/(137/60)/8}
\pgfmathsetmacro{\xC}{\xB+13*5/(137/60)/9}
\pgfmathsetmacro{\xD}{\xC+13*5/(137/60)/16}
\pgfmathsetmacro{\xE}{\xD+13*5/(137/60)/25}
\filldraw[zone,hatched] (\xA,0) rectangle (\xB,{2*\u});
\filldraw[zone,hatched] (\xA,{2*\u}) rectangle (\xB,{4*\u});
\filldraw[zone,fill=white] (\xA,{4*\u}) rectangle (\xB,{5*\u});
\node[lab] at ({(\xA+\xB)/2},{\u}) {$a=2$};
\node[lab] at ({(\xA+\xB)/2},{3*\u}) {$a=2$};
\filldraw[zone,hatched] (\xB,0) rectangle (\xC,{3*\u});
\filldraw[zone,fill=white] (\xB,{3*\u}) rectangle (\xC,{4*\u});
\filldraw[zone,fill=white] (\xB,{4*\u}) rectangle (\xC,{5*\u});
\node[lab] at ({(\xB+\xC)/2},{1.5*\u}) {$a=3$};
\filldraw[zone,hatched] (\xC,0) rectangle (\xD,{4*\u});
\filldraw[zone,fill=white] (\xC,{4*\u}) rectangle (\xD,{5*\u});
\node[lab] at ({(\xC+\xD)/2},{2*\u}) {$a=4$};
\filldraw[zone,hatched] (\xD,0) rectangle (\xE,{5*\u});
\node[lab,align=center] at ({(\xD+\xE)/2},{2.5*\u}) {$a$\\$=5$};
\foreach \j in {0,...,4}{
  \filldraw[zone,fill=white] (\xE,{\j*\u}) rectangle (\xF,{(\j+1)*\u});
}
\draw[line width=1pt] (\xA,0) rectangle (\xF,{5*\u});
\node[font=\small] at (-0.4,{2.5*\u}) {$m$};
\node[font=\small] at ({(\xA+\xF)/2},-0.4) {$n$};
\end{tikzpicture}
\caption{The patterns $Q_2,Q_3,Q_4,Q_5,Q_1$ of Example~\ref{ex:columns} for $m=5$, laid out as slabs of lengths proportional to the weights $\omega_a$, as in the proof of Theorem~\ref{thm:ratio}. Hatched slots receive the rectangles of height $a\ge2$. The white slots are unit slots.}
\label{fig:columns}
\end{figure}
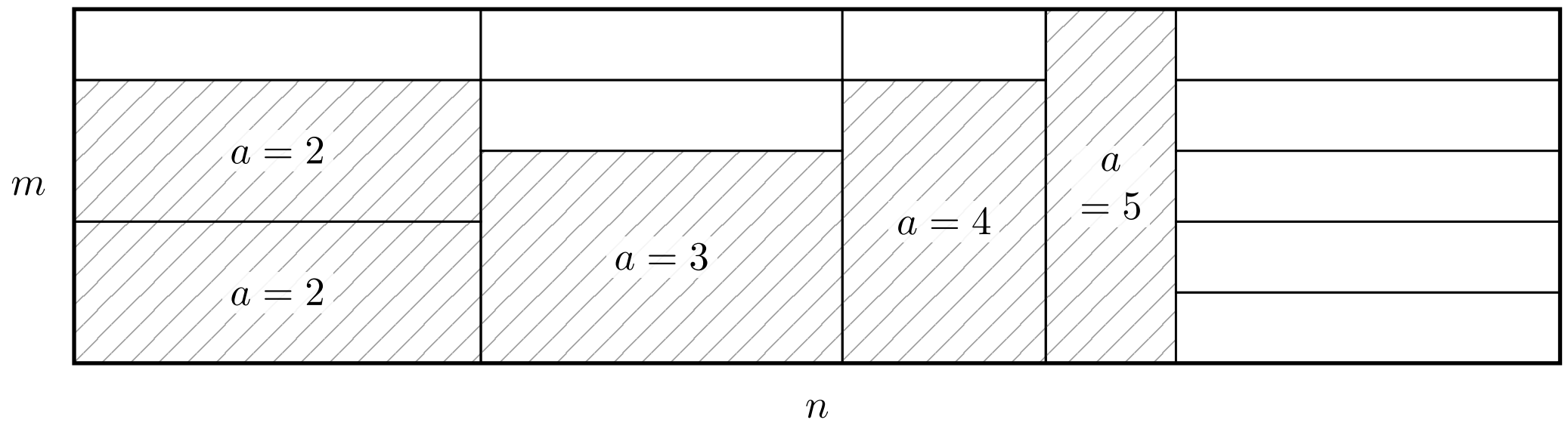

In particular, Lemma \ref{lm:patterns} ensures that there is an explicit choice for rectangles, shown in Figure~\ref{fig:columns}.

\begin{ex}\label{ex:columns}
Let $r=1$ and $\bh=(m)$ with $m\ge2$. For $2\le a\le m$, let $Q_a$ consist of $\lfloor m/a\rfloor$ slots of height $a$ followed by $m-a\lfloor m/a\rfloor$ unit slots, and let $Q_1$ consist of $m$ unit slots. Give $Q_a$ the weights
\[
\omega_a=\frac{m}{H_ma^2\lfloor m/a\rfloor}\quad(2\le a\le m)\qquad\text{and}\qquad\omega_1=1-\sum_{a=2}^m\omega_a.
\]
By Lemma~\ref{lm:columns} below, $\omega_1>1/(mH_m)>0$. For $a\ge2$ we get $\EE\nu_a=\omega_a\lfloor m/a\rfloor=m/(H_ma^2)=x_a$. Because every pattern has total height $m$, this also forces $\EE\nu_1=m-\sum_{a\ge2}ax_a=m/H_m=x_1$.
\end{ex}

\begin{lm}\label{lm:columns}
For every integer $m\ge2$,
\[
\sum_{a=2}^m\frac{m-a\lfloor m/a\rfloor}{a^2\lfloor m/a\rfloor}<1-\frac1m,
\qquad\text{equivalently}\qquad
\frac m{H_m}\sum_{a=2}^m\frac1{a^2\lfloor m/a\rfloor}<1-\frac1{mH_m}.
\]
\end{lm}
\begin{proof}
The equivalence follows from $m/(a^2\lfloor m/a\rfloor)=1/a+(m-a\lfloor m/a\rfloor)/(a^2\lfloor m/a\rfloor)$ and $\sum_{a=2}^m1/a=H_m-1$. Let $h=\lceil m/2\rceil$. For $2\le a\le h$ the remainder $m-a\lfloor m/a\rfloor$ lies in $\{0,\ldots,a-1\}$, and $z\mapsto z/(m-z)$ is increasing, so the summand is at most $(a-1)/(a(m-a+1))<1/(m-a+1)$. For $h<a\le m$ we have $\lfloor m/a\rfloor=1$ and $(m-a)/a^2<m/(a(a-1))-1/a$. The second range is nonempty, so
\[
\sum_{a=2}^m\frac{m-a\lfloor m/a\rfloor}{a^2\lfloor m/a\rfloor}<\bigl(H_{m-1}-H_{m-h}\bigr)+m\Bigl(\frac1h-\frac1m\Bigr)-\bigl(H_m-H_h\bigr)=H_h-H_{m-h}+\frac mh-1-\frac1m.
\]
Since $2h-m\in\{0,1\}$, we have $H_h-H_{m-h}=(2h-m)/h$, and the right-hand side equals $1-1/m$.
\end{proof}

\begin{re}
A Riemann-sum argument shows that the sum on the left of Lemma~\ref{lm:columns} tends to Euler's constant $\gamma$. Indeed, it equals $\frac1m\sum_{a\le m}f(a/m)$ with $f(x)=1/(x^2\lfloor1/x\rfloor)-1/x\in[0,2)$, and $\int_0^1f=\int_1^\infty(1/\lfloor y\rfloor-1/y)\,dy=\gamma$. Hence $\omega_1\sim(1-\gamma)/H_m$.
\end{re}

\section{Fixed cross-sections}\label{sec:fixed-proof}

This section is devoted to the proof of Theorem \ref{thm:fixed} and some consequences arising from it.
By Proposition~\ref{prop:coef} applied with $N=Vn$, it is enough to show the following.

\begin{thm}\label{thm:ratio}
As $n\to\infty$, we have $p(\bh,n)/U_{\bh}(Vn)\to1$.
\end{thm}

\begin{proof}
If $\bh=(1,\ldots,1)$, every multiset of volume $n$ is a partition of $n$ and tiles the host, so $p(\bh,n)=U_{\bh}(n)$. Assume from now on that $h_r\ge2$.

Take the random pattern $\mathsf P$ of Lemma~\ref{lm:patterns}, or of Example~\ref{ex:columns} when $r=1$. Let $\omega_0=\PP(\mathsf P\text{ is all-unit})$. Then $0<\omega_0<1$, since $\EE\nu_a(\mathsf P)=x_a>0$ for $a\ne\e$. Let $\cP$ be the finite set of patterns other than the all-unit one that occur with positive probability, and put $\lambda_P=\PP(\mathsf P=P)/(1-\omega_0)$ for $P\in\cP$. For $a\ne\e$, the all-unit pattern has $\nu_a=0$, so
\begin{equation}\label{eq:excess}
\sum_{P\in\cP}\lambda_P\nu_a(P)=\frac{x_a}{1-\omega_0}.
\end{equation}
Fix $\varepsilon>0$ with $2\varepsilon<\omega_0x_a/(1-\omega_0)$ for all $a\ne\e$, and put $R=V(|\cP|+1)$.

Let $N=Vn$, $t=t_N$, and let $L$ and $w$ be as in Lemma~\ref{lm:typical}. Let $\mathcal G_n$ be the set of multisets of volume $N$ satisfying events (i)--(iv) of that lemma with these $\varepsilon$ and $R$. By Lemma~\ref{lm:typical} and uniformity, $|\mathcal G_n|=(1-o(1))\,U_{\bh}(N)$. Since $L=O(\sqrt n\log n)$, every piece of a multiset in $\mathcal G_n$ fits individually once $n$ is large. It therefore suffices to show that, for all large $n$, every $\mu\in\mathcal G_n$ tiles the host $[0,h_1]\times\cdots\times[0,h_r]\times[0,n]$.

\emph{Slabs and tubes.} Divide the axial interval $[0,n]$ into consecutive intervals of lengths $\lfloor\lambda_Pn\rfloor$ for $P\in\cP$, followed by a remainder of length less than $|\cP|$. Use pattern $P$ throughout its interval, so each slot of $P$ becomes a \emph{tube}: the slot times the interval. The tubes of footprint $a\ne\e$ have total length
\[
\ell_a=\sum_{P\in\cP}\nu_a(P)\lfloor\lambda_Pn\rfloor\ge\frac{x_an}{1-\omega_0}-O(1)\ge(x_a+2\varepsilon)n-O(1)
\]
by~\eqref{eq:excess}, and their number is bounded independently of $n$.

\emph{Non-unit footprints.} Fix $a\ne\e$. The pieces of footprint $a$ have lengths at most $L$, by (i), and total length $W_a\le(x_a+\varepsilon)n$, by (ii). Since $\ell_a-W_a\ge\varepsilon n-O(1)$ is much larger than $L$ times the number of tubes, Lemma~\ref{lm:nextfit} places them end to end in the tubes of footprint $a$, with sides $a_1,\ldots,a_r$ along the axes $1,\ldots,r$.

\emph{Unit footprint.} Split the unused end of every tube of footprint $a\ne\e$ into $v_a$ parallel tubes of unit footprint. Together with the original unit tubes and the $V$ unit tubes over the remainder interval, this leaves at most $\sum_{P\in\cP}\sum_av_a\nu_a(P)+V=R$ tubes of unit footprint (we add tubes of length $0$ if necessary). Their total length is the volume not yet covered, namely $N-\sum_{a\ne\e}v_aW_a=W_{\e}$. The pieces of footprint $\e$ are rods $1\times\cdots\times1\times c$. Their lengths form a multiset with sum $W_\e$ and all elements at most $L$, by (i). The elements in $\{2,\ldots,w\}$ have sum greater than $RL$, by (iii), and at least $Rw$ elements equal $1$, by (iv). Lemma~\ref{lm:exact}, with $k=R$, splits them exactly among the unit tubes.

All coordinates are integers, so $\mu$ tiles the host. Hence $p(\bh,n)\ge(1-o(1))U_{\bh}(Vn)$, and with $p(\bh,n)\le U_{\bh}(Vn)$ this proves Theorem~\ref{thm:ratio}.
\end{proof}

As an illustration of Theorem \ref{thm:fixed}, let us consider the box partition function $p(2,2,n)$.

\begin{ex}
For $\bh=(2,2)$ the footprints are $(1,1),(1,2),(2,2)$, with $v_a=1,2,4$ and $b_a=1,2,2$. Hence $\Hh=7/4$, $B=7/2$ and $C=8/\pi^{3/2}$, and Theorem~\ref{thm:fixed} gives
\begin{equation}\label{eq:22}
p(2,2,n)\sim\frac{49\pi^2}{4608\,n^{5/2}}\exp\left(\pi\sqrt{\frac{14n}{3}}\right).
\end{equation}
\end{ex}

\subsection{Proof of Theorem~\ref{cor:rect}}
Let $\bh=(m)$. Then $\mathcal{A}_{\bh}=\{1,\ldots,m\}$ with $v_a=b_a=a$, the type $(a,c)$ is the rectangle $a\times c$, and $U_{\bh}(N)=U_m(N)$, where
\begin{equation}\label{eq:Um}
U_m(N)=[x^N]\prod_{a=1}^{m}\prod_{b=a}^{\infty}\frac1{1-x^{ab}}
\end{equation}
is the number of multisets of rectangles $a\times b$ with $a\le m$ and total area $N$. Further $\Hh=H_m$, $A=\pi^2H_m/6$, $B=\sum_{a=1}^m(a-\frac12)=m^2/2$ and
\[
C=\prod_{a=1}^m\frac{a^{a-\tfrac12}(a-1)!}{\sqrt{2\pi}}=\frac{(m!)^{m-\tfrac12}}{(2\pi)^{m/2}},
\]
because $\prod_{a=1}^m(a-1)!=\prod_{a=1}^ma^{m-a}$. Theorem~\ref{thm:fixed} with $N=mn$ gives
\[
p(m,n)\sim U_m(mn)\sim\frac{(m!)^{m-\tfrac12}}{2\sqrt\pi(2\pi)^{m/2}}\left(\frac{\pi^2H_m}6\right)^{\frac{m^2+1}4}(mn)^{-\frac{m^2+3}4}\exp\left(\pi\sqrt{\frac{2H_mmn}3}\right),
\]
which is Theorem~\ref{cor:rect}. For $m=1,2,3$ one has $H_m=1,\frac32,\frac{11}6$, and substitution gives the three special cases discussed after it.\qed

\subsection{Multiplicities of fixed types}
The proof of Theorem \ref{thm:ratio} also describes a uniformly random partition. For a type $\theta=(a,c)$ write $|\theta|=v_ac$ for its volume. Let $\PP_{\rm tile}$ be the uniform distribution on the $p(\bh,n)$ partitions of the host, let $\PP_{\rm all}$ be the uniform distribution on the $U_{\bh}(Vn)$ multisets of volume $Vn$, and let $X_\theta$ be the multiplicity of $\theta$.

\begin{cor}\label{cor:multiplicities}
Fix $\bh$ and distinct types $\theta_1,\ldots,\theta_k$, and let $t_N=\sqrt{A/N}$ with $N=Vn$. Under $\PP_{\rm tile}$, as $n\to\infty$,
\[
\bigl(t_N|\theta_1|X_{\theta_1},\ldots,t_N|\theta_k|X_{\theta_k}\bigr)\xrightarrow{\ d\ }(E_1,\ldots,E_k),
\]
where $E_1,\ldots,E_k$ are independent exponential random variables of mean $1$. For rectangles, $t_N|\theta|=\pi ab\sqrt{H_m/(6mn)}$ for the type $a\times b$.
\end{cor}

\begin{proof}
For any event $E$, decomposing $\PP_{\rm all}(E)$ according to tileability gives
\[
|\PP_{\rm tile}(E)-\PP_{\rm all}(E)|\le\PP_{\rm all}(\text{not tileable})=1-\frac{p(\bh,n)}{U_{\bh}(Vn)},
\]
which tends to zero by Theorem~\ref{thm:ratio}. So it suffices to work under $\PP_{\rm all}$. Fix $y_1,\ldots,y_k\ge0$ and put $\kappa_i=\lfloor y_i/(t_N|\theta_i|)\rfloor+1$, so that $\{t_N|\theta_i|X_{\theta_i}>y_i\}=\{X_{\theta_i}\ge\kappa_i\}$. For $I\subseteq\{1,\ldots,k\}$, deleting $\kappa_i$ copies of $\theta_i$ for $i\in I$ is a bijection onto the multisets of volume $N-J_I$, where $J_I=\sum_{i\in I}\kappa_i|\theta_i|$. Hence
\[
\PP_{\rm all}\bigl(X_{\theta_i}\ge\kappa_i\ \forall i\in I\bigr)=\frac{U_{\bh}(N-J_I)}{U_{\bh}(N)}.
\]
We have $t_NJ_I\to\sum_{i\in I}y_i$, so $J_I=O(\sqrt N)$, and Proposition~\ref{prop:coef} gives
\[
\frac{U_{\bh}(N-J_I)}{U_{\bh}(N)}=(1+o(1))\exp\left(-\frac{2\sqrt AJ_I}{\sqrt N+\sqrt{N-J_I}}\right)\longrightarrow\exp\Bigl(-\sum_{i\in I}y_i\Bigr).
\]
Inclusion--exclusion over $I$ gives $\PP_{\rm all}(t_N|\theta_i|X_{\theta_i}\le y_i\ \forall i)\to\prod_i(1-e^{-y_i})$, which is the claim.
\end{proof}

\begin{re}
This is an instance of the independent-process approximations of Arratia and Tavar\'e~\cite{AT94}. The Boltzmann model of Section~\ref{sec:boltzmann} provides the independent process, and Theorem~\ref{thm:ratio} transfers the conclusion from all multisets to tiling multisets.
\end{re}

\section{Several growing sides}\label{sec:growing}
Our main goal is to exhibit the proof of Theorem~\ref{thm:growing}. Throughout this section we use the notation from its statement. 
When $r=0$, we put $\mathcal{A}_{\bh}=\{\varnothing\}$,
$\e=\varnothing$ and $v_\varnothing=b_\varnothing=1$,
so that $V=\Hh=1$; for $r\ge1$, we write
$V=h_1\cdots h_r$, $v_a=a_1\cdots a_r$ and $b_a=a_r$
as before.
Put $H_i=\floor{n^{\alpha_i}}$, $G=H_1\cdots H_g$
and $\beta=\frac12\sum_{i=1}^g\alpha_i$, so that $N=VG$.
If $r\ge1$, let $n$ be large enough that $h_r\le H_1$.
By Lemma~\ref{lm:fit}, the pieces that fit the host
individually are the boxes of type $(a,j)$,
where $a\in\mathcal{A}_{\bh}$ and $j=(j_1,\ldots,j_g)$
satisfies $b_a\le j_1\le\cdots\le j_g$ and
$j_i\le H_i$ for all $i\in\{1,\ldots,g\}$;
the box has sorted sides $(a_1,\ldots,a_r,j_1,\ldots,j_g)$
and volume $v_aj_1\cdots j_g$. Let us put
\[
\tau=\sqrt{\frac{\z\,\Hh\,D_\alpha(\log n)^{g-1}}N},\qquad J=N\tau=\sqrt{\z\,\Hh\,D_\alpha N}\,(\log n)^{\frac{g-1}2},\qquad\rho_a=\frac1{\Hh v_a},
\]
so that $\sum_a\rho_a=1$. Since $2J=2\sqrt{\z\Hh D_\alpha N}(\log n)^{(g-1)/2}=\pi\sqrt{2\Hh D_\alpha N/3}\,(\log n)^{(g-1)/2}$, Theorem~\ref{thm:growing} asserts $\log p\sim2J$. The positivity of $D_\alpha$ follows because the point $y_i=\alpha_i/2+\eta\,(i-\frac{g+1}2)$ lies in the interior of $\mathcal D_\alpha$ for small $\eta>0$. Note that $\log\tau/\log n\to-\beta$.

For $\kappa\in\{0,g+2\}$ put $H_i^{(\kappa)}=\lfloor H_i/(\log n)^\kappa\rfloor$, and for $t>0$ let
\[
Z_a^{(\kappa)}(t)=\prod_{\substack{b_a\le j_1\le\cdots\le j_g\\ j_i\le H_i^{(\kappa)}}}\bigl(1-e^{-tv_aj_1\cdots j_g}\bigr)^{-1}.
\]
The full caps $\kappa=0$ describe all individually fitting pieces. The reduced caps $\kappa=g+2$ are used for the lower bound.

\subsection{The generating product}
In order to prove Theorem \ref{thm:growing}, we need the forthcoming auxiliary result, which plays an important role in both the upper and lower estimates.
\begin{lm}\label{lm:growing-product}
For fixed $a\in\mathcal{A}_{\bh}$, fixed $s>0$ and $\kappa\in\{0,g+2\}$,
\[
\log Z_a^{(\kappa)}(s\tau)=(1+o(1))\,\frac{\rho_a}{s}\,J\qquad(n\to\infty).
\]
The same holds if the single factor of the unit type ($a=\e$, $j=(1,\ldots,1)$) is omitted.
\end{lm}

\begin{proof}
Let $t=s\tau$ and $\phi(x)=-\log(1-e^{-x})$. This function is positive and decreasing on $(0,\infty)$, with $\int_0^\infty\phi=\z$. Write $j'=(j_1,\ldots,j_{g-1})$, $\Pi=j_1\cdots j_{g-1}$ and $\vartheta=tv_a\Pi$. Summing over $j_g$ first gives
\begin{equation}\label{eq:riemann}
tv_a\log Z_a^{(\kappa)}(t)=\sum_{j'}\frac1{\Pi}\,\Phi(j'),\qquad\Phi(j')=\vartheta\sum_{k=j_{g-1}}^{H_g^{(\kappa)}}\phi(\vartheta k),
\end{equation}
where $j'$ runs over $b_a\le j_1\le\cdots\le j_{g-1}$ with $j_i\le H_i^{(\kappa)}$ for $i<g$ (and $j_{g-1}\le H^{(\kappa)}_g$). Since $\phi$ is decreasing,
\begin{equation}\label{eq:bracket}
\int_{\vartheta j_{g-1}}^{\vartheta(H_g^{(\kappa)}+1)}\phi\;\le\;\Phi(j')\;\le\;\min\Bigl\{\z,\ \vartheta\phi(\vartheta j_{g-1})+\int_{\vartheta j_{g-1}}^\infty\phi,\ \int_0^{\vartheta H_g^{(\kappa)}}\phi\Bigr\}.
\end{equation}

Use the logarithmic coordinates $y_i=\log j_i/\log n$ for $i<g$, and set $y_g=\beta-\sum_{i<g}y_i$. Since $\log t=-\beta\log n+O(\log\log n)$ and $\log H_g^{(\kappa)}=\alpha_g\log n+O(\log\log n)$, we have, uniformly in $j'$,
\[
\vartheta j_{g-1}=n^{\,y_{g-1}-y_g+o(1)},\qquad\vartheta H_g^{(\kappa)}=n^{\,\alpha_g-y_g+o(1)}.
\]
Fix $\delta>0$. By~\eqref{eq:bracket}, $\Phi\to\z$ uniformly on $\{y_{g-1}\le y_g-\delta,\ y_g\le\alpha_g-\delta\}$. Also $\Phi\to0$ uniformly on $\{y_{g-1}\ge y_g+\delta\}$, using the middle term of the minimum, $\vartheta\le\vartheta j_{g-1}$ and $x\phi(x)\to0$; and $\Phi\to0$ uniformly on $\{y_g\ge\alpha_g+\delta\}$, using the last term.

For fixed $0\le u<v$ we have $\sum_{n^u\le j<n^v}1/j=(v-u)\log n+O(1)$. Hence, for every fixed box $Q\subset[0,\infty)^{g-1}$, the normalised harmonic weight $(\log n)^{-(g-1)}\sum_{j':\,y\in Q}\Pi^{-1}$ tends to the volume of $Q$, where the sum is taken over all positive integer tuples $j'=(j_1,\ldots,j_{g-1})$ with $y\in Q$. All relevant $y$ lie in the bounded set $\prod_{i<g}[0,\alpha_i+1]$. Partition this set into cubes of side $\delta$. Consider a cube at distance more than $\delta$ from all the hyperplanes
\[
y_1=0,\quad y_i=y_{i+1}\ (i<g-1),\quad y_i=\alpha_i\ (i<g),\quad y_{g-1}=y_g,\quad y_g=\alpha_g.
\]
For large $n$, such a cube lies entirely inside or entirely outside the range of summation, and $\Phi$ is uniformly within $o(1)$ of $\z$ (if the cube lies in $\mathcal D_\alpha$) or of $0$ (otherwise). The remaining cubes have total volume $O(\delta)$, and there $0\le\Phi\le\z$. The terms with $j_1<b_a$ lie in the neighbourhood of $y_1=0$ and are covered by the same bound. Therefore
\[
(\log n)^{-(g-1)}\,tv_a\log Z_a^{(\kappa)}(t)=\z D_\alpha+O(\delta)+o(1).
\]
As $\delta>0$ is arbitrary,
\[
\log Z_a^{(\kappa)}(s\tau)=(1+o(1))\,\frac{\z D_\alpha(\log n)^{g-1}}{s\tau v_a}.
\]
Since $\z D_\alpha(\log n)^{g-1}=\tau^2N/\Hh$, this equals $(1+o(1))J\rho_a/s$. The unit type contributes $\phi(t)=O(\log n)=o(J)$.
\end{proof}

\subsection{Upper bound}
Let $\widetilde U(N)=[q^N]\prod_aZ_a^{(0)}(q)$ denote the number of multisets of volume $N$ of individually fitting pieces, where $Z_a^{(0)}(q)$ denotes the same product as a power series in $q$. Then $p\le\widetilde U(N)$, and positivity of coefficients together with Lemma~\ref{lm:growing-product} gives
\[
\log\widetilde U(N)\le\tau N+\sum_a\log Z_a^{(0)}(\tau)=J+(1+o(1))J\sum_a\rho_a=(2+o(1))J.
\]

\subsection{Lower bound}
Fix $0<u<1$, put $c=u^{-\tfrac12}$ and $t=c\tau$. Consider the types with reduced caps $\kappa=g+2$, excluding the unit type. Give them independent geometric multiplicities $X_{a,j}$ with parameters $e^{-tv_aj_1\cdots j_g}$, as in Section~\ref{sec:boltzmann}, and let $Z(t)=\prod_aZ_a^{(g+2)}(t)$ with the unit factor omitted. A multiset $\mu$ then has probability $e^{-tT(\mu)}/Z(t)$, where $T(\mu)$ is its volume. Let $T_a$ be the volume of the pieces with footprint $a$, and let $Y_a=T_a/v_a$ be their volume in the growing coordinates.
Let $Z_a(t)$ denote $Z_a^{(g+2)}(t)$ with
the unit factor omitted when $a=\e$.\\
\emph{Concentration.} For fixed $|h|<c$, independence and Lemma~\ref{lm:growing-product} give
\[
\frac1J\log\EE_te^{h\tau T_a}=\frac1J\log\frac{Z_a((c-h)\tau)}{Z_a(c\tau)}\longrightarrow\Lambda_a(h):=\rho_a\Bigl(\frac1{c-h}-\frac1c\Bigr).
\]
Here $\Lambda_a(0)=0$ and $\Lambda_a'(0)=\rho_a/c^2=u\rho_a$. Given $\varepsilon'>0$, choose a small $h>0$ with $\Lambda_a(h)\le h(u\rho_a+\varepsilon'/2)$ and $\Lambda_a(-h)\le-h(u\rho_a-\varepsilon'/2)$. Markov's inequality then gives
\[
\PP_t\bigl(T_a\ge(u\rho_a+\varepsilon')N\bigr)\le e^{-h\tau(u\rho_a+\varepsilon')N}\EE_te^{h\tau T_a}\le e^{-(h\varepsilon'/2+o(1))J}\to0,
\]
and the lower tail is handled in the same way. Put $x_a=\rho_aV/v_a=V/(S_{\bh}v_a^2)$.
For every $\varepsilon>0$ the event
\[
\mathcal E=\Bigl\{T\ge(u-\varepsilon)N\ \text{ and }\ Y_a\le(u+\varepsilon)x_aG\ \text{for every }a\Bigr\}
\]
has probability tending to one. Every multiset in $\mathcal E$ has probability at most $e^{-t(u-\varepsilon)N}/Z(t)$, so Lemma~\ref{lm:growing-product} gives
\begin{equation}\label{eq:count-E}
\log|\mathcal E|\ge\log Z(c\tau)+c\tau(u-\varepsilon)N+o(1)=\Bigl(\frac1c+c(u-\varepsilon)+o(1)\Bigr)J.
\end{equation}

\emph{Bins.} If $r\ge1$, let $\cP$ be the support of the pattern law of Lemma~\ref{lm:patterns}, with weights $\lambda_P=\PP(\mathsf P=P)$; if $r=0$, let $\cP$ consist of the empty pattern with one slot of footprint $\varnothing$ and $\lambda=1$, so that $x_\varnothing=1$. Divide the last growing axis $[0,H_g]$ into consecutive intervals of lengths $\lfloor\lambda_PH_g\rfloor$, and use pattern $P$ on the fixed coordinates over its interval. Each slot of footprint $a$ then gives a \emph{bin}: the slot, times $[0,H_1]\times\cdots\times[0,H_{g-1}]$, times the interval. The bins of footprint $a$ have total height
\[
\sum_P\nu_a(P)\lfloor\lambda_PH_g\rfloor\ge x_aH_g-O(1)
\]
along axis $g$, and their number is bounded.

\emph{Layers.} Fix $\delta>0$ with $(1+\delta)^g(u+\varepsilon)<1$, and let $\mu\in\mathcal E$. Group the pieces of $\mu$ with footprint $a$ into classes according to the vector $(k_1,\ldots,k_g)$ with $j_i\in[(1+\delta)^{k_i},(1+\delta)^{k_i+1})$. For a class, let $d_i$ be the largest $j_i$ occurring in it, so that $d_i<(1+\delta)j_i$ for every piece of the class, and $d_i\le H_i^{(g+2)}\le H_i/(\log n)^{g+2}$. Since every $j_i\le n^{\alpha_g}$, there are $O_\delta((\log n)^g)$ classes.

A \emph{layer} for a class is the slot, times a grid of $Q=\prod_{i<g}\lfloor H_i/d_i\rfloor$ cells of size $d_1\times\cdots\times d_{g-1}$ in $[0,H_1]\times\cdots\times[0,H_{g-1}]$, times an interval of length $d_g$. Each cell receives one piece, with its footprint sides along the fixed axes and $j_i\le d_i$ along growing axis $i$. A class with $q$ pieces uses $\lceil q/Q\rceil$ layers. Put $G'=H_1\cdots H_{g-1}$. Since $d_i/H_i\to0$ uniformly, $Q\ge(1-o(1))G'/(d_1\cdots d_{g-1})$, and the total height of all layers of footprint $a$ is at most
\[
\sum_{\text{classes}}\Bigl(\frac qQ+1\Bigr)d_g\le(1+o(1))(1+\delta)^g\frac{Y_a}{G'}+O_\delta\Bigl(\frac{H_g}{(\log n)^2}\Bigr)\le\bigl((1+\delta)^g(u+\varepsilon)+o(1)\bigr)x_aH_g.
\]
This is less than the total bin height $x_aH_g-O(1)$ by a positive multiple of $H_g$. Each layer has height at most $\max d_g=o(H_g)$, and the number of bins is bounded. Lemma~\ref{lm:nextfit}, applied with heights in place of lengths, therefore stacks all layers into the bins of footprint $a$.

\emph{Injectivity.} Fill every uncovered unit cell of the host with a unit cube. The result is a partition of the host consisting of $\mu$ and some unit cubes. Since $\mu$ contains no unit cubes, removing them recovers $\mu$. Thus $p\ge|\mathcal E|$, and by~\eqref{eq:count-E}
\[
\liminf_{n\to\infty}\frac{\log p}{J}\ge\frac1c+c(u-\varepsilon)=2\sqrt u-\frac{\varepsilon}{\sqrt u}.
\]
Letting $\varepsilon\downarrow0$ and then $u\uparrow1$ gives $\liminf\log p/J\ge2$. Together with the upper bound, and with $p\le\widetilde U(N)$, this proves Theorem~\ref{thm:growing}.\qed

\subsection{Proofs of Corollaries~\ref{cor:cubes} and~\ref{cor:special}, and further examples}
Now, we exhibit some applications of Theorem~\ref{thm:growing}. In particular, we prove both Corollary~\ref{cor:cubes} and Corollary~\ref{cor:special}. 

\begin{proof}[Proof of Corollary~\ref{cor:cubes}]
Here $r=0$ and $\alpha=(1,\ldots,1)$, so $\beta=d/2$ and $N=n^d$. The set $\{y\in[0,1]^d:\sum y_i=d/2\}$, projected to its first $d-1$ coordinates, has $(d-1)$-volume equal to the density at $d/2$ of a sum of $d$ independent uniform $[0,1]$ variables, namely
\[
\frac1{(d-1)!}\sum_{j\ge0}(-1)^j\binom dj\Bigl(\frac d2-j\Bigr)_+^{d-1}.
\]
Coordinate permutations preserve this projected measure (it is the Euclidean surface measure divided by $\sqrt d$), and the sets where two coordinates coincide have measure zero. So the sorted region $\mathcal D_\alpha$ has volume
\[
D_\alpha=\frac{1}{d!\,(d-1)!}\sum_{j=0}^{\lfloor d/2\rfloor}(-1)^j\binom dj\Bigl(\frac d2-j\Bigr)^{d-1},
\]
and Theorem~\ref{thm:growing} gives $\log p\sim\pi\sqrt{2D_\alpha/3}\;n^{d/2}(\log n)^{(d-1)/2}$, that is, the constant $\kappa_d$. For $d=2,3,4$ one finds $D_\alpha=\frac12,\frac18,\frac1{36}$, hence $\kappa_2=\pi/\sqrt3$, $\kappa_3=\pi/(2\sqrt3)$ and $\kappa_4=\pi/(3\sqrt6)$, as displayed before Corollary~\ref{cor:cubes}. For $d=5$ one finds $D_\alpha=\frac{23}{4608}$ and $\kappa_5=\pi\sqrt{23}/(48\sqrt3)$.
\end{proof}

For two growing sides, $g=2$ and $\alpha=(s,u)$ with $s\le u$, we have $\beta=(s+u)/2$, and $\mathcal D_\alpha$ is the set of $y_1$ with $0\le y_1\le\beta-y_1$, $y_1\le s$ and $\beta-y_1\le u$. The last condition reads $y_1\ge(s-u)/2$, which is automatic. Hence
\begin{equation}\label{eq:D2}
D_{(s,u)}=\min\Bigl\{s,\frac{s+u}4\Bigr\}.
\end{equation}

\begin{proof}[Proof of Corollary~\ref{cor:special}]
Here $r=0$, $N=hw$, and $\log h=(s+o(1))\log n$, $\log(hw)=(s+u+o(1))\log n$. If $u\le3s$ then $D_\alpha=(s+u)/4$ by~\eqref{eq:D2}, and
\[
2\sqrt{\z\,\tfrac{s+u}4\,hw\log n}=\sqrt{\z\,hw\log(hw)}\,(1+o(1)).
\]
If $u\ge3s$ then $D_\alpha=s$, and
\[
2\sqrt{\z\,s\,hw\log n}=2\sqrt{\z\,hw\log h}\,(1+o(1)).
\] Since $\sqrt{\z}=\pi/\sqrt6$ and $2\sqrt{\z}=\pi\sqrt{2/3}$, this is the claim.
\end{proof}

\begin{ex}\label{ex:more}
(a) For fixed $k\ge1$, take $r=1$, $\bh=(k)$ and $\alpha=(1,1)$. Then $V=k$, $S_\bh=H_k$ and $D_\alpha=\frac12$ by~\eqref{eq:D2}, so
\[
\log p(k,n,n)\sim2\sqrt{\z H_k\cdot\tfrac12\cdot kn^2\log n}=\pi\sqrt{\frac{kH_k}{3}}\;n\sqrt{\log n}.
\]
For $k=1$ this is the case $d=2$ of Corollary~\ref{cor:cubes}, since $p(1,n,n)=p(n,n)$.

(b) Take $r=0$ and $\alpha=(1,1,2)$. Then $\beta=2$ and $\mathcal D_\alpha=\{0\le y_1\le y_2,\ y_1+2y_2\le2\}$; the conditions $y_2\le1$ and $y_3\le2$ are implied. This triangle has vertices $(0,0)$, $(0,1)$, $(\frac23,\frac23)$ and area $\frac13$. With $N=n^4$ we get
\[
\log p(n,n,n^2)\sim2\sqrt{\z n^4/3}\,\log n=\frac{\pi\sqrt2}3\,n^2\log n.
\]
\end{ex}

\section{Further problems}\label{sec:problems}

In the hopes of encouraging further research on this topic, we pose the following further problems:  

\begin{enumerate}
\item \emph{Error terms.} Determine the rate in $p(m,n)/U_m(mn)\to1$ and an effective error term in Theorem~\ref{cor:rect}. 
Determine also the range of $m=m(n)$ for which Theorem~\ref{cor:rect} remains valid.
\item \emph{Squares and cubes.} Find lower-order terms and multiplicative asymptotics for $p(n,n)$ and $p(n,n,n)$. The data of~\cite[Table~2]{GZsquare} suggest that $p(n,n)$ is asymptotic to the corresponding relaxed count considered there; our results only show that the two agree to leading order in the logarithm. The analogous question for the ratio
\[
\frac{p(n,n,n)}{[x^{n^3}]\prod_{1\le a\le b\le c\le n}(1-x^{abc})^{-1}}
\]
is open. Our method does not resolve the limiting behaviour of this ratio.
At 
$t\asymp(\log n)n^{-3/2}$ there are $\asymp1/t$ types $1\times b\times c$ with $b\asymp1/(tn)$ and $c\asymp n$, each present with probability bounded below. So a typical multiset contains many pieces of macroscopic length, which the tail bounds of Section~\ref{sec:fixed-proof} unfortunately do not control.
\end{enumerate}

\section*{Statement on AI usage}

The initial ideas for the main results and their proofs are due to the authors. This paper was prepared with the assistance of ChatGPT 6.2 (OpenAI) and Claude Opus 5.5 (Anthropic).

\bibliographystyle{abbrv}  
\bibliography{bib}

\end{document}